\documentclass{amsart}
\usepackage{amsfonts}
\usepackage{amsmath}
\usepackage{amscd}
\usepackage{amssymb}

\def\?{\char'76}
\def\!{\char'74}

\def\n{\mathfrak{n}}

\newcounter{numero}

\newcounter{letra}

\theoremstyle{definition}

\theoremstyle{remark}

\begin{document}

\title{On Normal Stratified Pseudomanifolds}

\author{G. Padilla}

\date{May 9/2002, revised on June 15/2002}

\address{Universidad Central de Venezuela, Caracas}
\email{gabrielp@euler.ciens.ucv.ve}

\thanks{This paper was supported by the the ECOSNord project,
Euskal Herria Unibertsitatea's Math Department, the University of
Artois' LaboGA and the Universidad Central de Venezuela's C.D.C.H.
The author would like to thank M. Saralegi for encouraging us to
look at this problem, and also R. Popper for some helpful
conversations.}



\dedicatory{Devoted to the victims of the natural
tragedy in Vargas, Dec.15/1999,\\
who died under the rage of Waraira Repano.}

\begin{abstract}
    Any $pl$-stratified pseudomanifod can be normalized
    preserving its intersection homology. In this paper  we extend
    this result for any topological stratified pseudomanifold and for
    a family of perversities which is larger than usual. Our
    construction is functorial. We also give a detailed description of
    the normalizer's stratification in terms of the initial stratified
    pseudomanifold.
\end{abstract}

\maketitle

\section*{Foreword}

{\it For an entire version of this article the reader should go better to
\tt Extracta Math. Vol. 15 (3), 383-412.}\vskip5mm

A stratified pseudomanifold is normal if its links are connected.
A normalization of a stratified pseudomanifold $X$ is a normal
stratified pseudomanifold $X^N$ together with a finite-to-one
projection $\n:X^N\rightarrow X$ preserving the intersection
homology. Recall that intersection homology is the suitable
algebraic tool for the stratified point of view: it was first
introduced by Goresky and MacPherson in the $pl$-category and
later extended for any topological stratified pseudomanifold
\cite{gm1}, \cite{gm2}. Following Borel the map $\n$ is usually
required to satisfy the following property: For each $x\in X$
there is a distinguished neighborhood $U$ such that the points of
$\n^{-1}(x)$ are in correspondence with the connected components
of the regular part of $U$. A normalization satisfying the above
condition always exists for any $pl$-stratified pseudomanifold
\cite{borel}, \cite{maccrory}. In this article we study the main
properties of the map $\n$. More precisely, we prove that $\n$ can
be required to satisfy a stronger condition: it is a locally
trivial stratified morphism preserving the conical structure
transverse to the strata. We make an explicit construction of such
a normalization for any topological stratified pseudomanifold. Our
construction is functorial, thus unique. We exhibit the
relationship between  the stratifications of $X$ and $X^N$.
Finally we prove that the normalization preserves the intersection
homology with the family of perversities given in
\cite{pervsheaf}, see also \cite{king},\cite{illinois}. This
family of perversities is larger than the usual one. With little
adjust our procedure holds also in the $C^{\infty}$ category.

\end{document}